\documentclass{article}[10pt]
\usepackage[a4paper, total={6.5in, 9in}]{geometry}
\usepackage{amsmath} 
\usepackage{amsthm}
\usepackage{graphicx} 
\usepackage{float} 
\usepackage{tikz}
\usepackage{multirow}
\usepackage{comment}

\newtheorem{theorem}{Theorem}
\newcommand{\tri}{\mathbin{\overline{\nabla}}}

\newcommand{\m}{{\mathbf m}}
\newcommand{\C}{{\mathbf c}}
\newcommand{\T}{{\mathbf t}}

\newcommand{\Sb}{{\mathbf S}}

\newcommand{\Sm}{\mathcal{S}}
\newcommand{\OB}{\big[}
\newcommand{\CB}{\big]}
\newcommand{\mC}{C_{\m}}
\newcommand{\mV}{V_{\m}}
\newcommand{\mE}{E_{\m}}
\newcommand{\mF}{F_{\m}}

\providecommand{\U}[1]{\protect\rule{.1in}{.1in}}

\title{\textbf{Finite Interpretations of a Hyper-Catalan Series Solution to Polynomial Equations and Visualizations}}
\author{Pratham Mukewar}
\date{March 20th, 2025}

\begin{document}

\maketitle

\begin{abstract}
The solution to the general univariate polynomial equation has been sought for centuries. It is well known there is no general solution in radicals for degrees five and above.

The hyper-Catalan numbers $C[m_2,m_3,m_4,\ldots]$ count the ways to subdivide a planar polygon into exactly $m_2$ triangles, $m_3$ quadrilaterals, $m_4$ pentagons, etc. Wildberger and Rubine (2025) show the generating series S of the hyper-Catalan numbers is a formal series zero of the general geometric polynomial (meaning, general except for a constant of $1$ and a linear coefficient of $-1$).

Using a variant of the series solution to the geometric polynomial that has the number of vertices, edges, and faces explicitly shown, We prove their infinite series result may be viewed as a finite identity at each level, where a level is a truncation of $\Sb$ to a given maximum number of vertices, edges, or faces (bounded by degree).

We illustrate this result, as well as the general correspondence between operations on sets of subdivided polygons and the algebra of polynomials, with figures and animations generated using Python.
   
\end{abstract}

\section{\fontsize{14}{16}\selectfont \textbf{Wildberger and Rubine's Polynomial Formula}}

In a paper published the \textit{American Mathematical Monthly}, Wildberger and Rubine \cite{Wildberger2025} (WR) solve the general geometric polynomial:

\begin{theorem}[The geometric polynomial formula] \label{thm:geompoly} The equation
$$ \displaystyle
0=1 - \alpha + t_2  \alpha ^2 + t_3  \alpha^3 + t_4 \alpha^4 + t_5 \alpha^5 \ + \  \ldots 
$$
has a formal power series solution:

\begin{align*}
\alpha &= 
\!\!
\sum_{\substack{m_k \ge 0}}  
\!
\dfrac{( 2m_2 + 3m_3 + 4m_4 + \ldots )! }{(1 + m_2 + 2m_3  +\ldots)!\, m_2! \, m_3! \cdots}  \,   t_2^{m_2 } t_3^{m_3}\cdots 
\end{align*}
\end{theorem}

WR call the coefficient the hyper-Catalan number $C[m_2, m_3, \ldots]$.
WR's result is that the generating function for the hyper-Catalan numbers is the solution to the general geometric polynomial.  Let's unwind what that means.

By \textbf{geometric polynomial}, they mean a univariate polynomial with a constant coefficient of $1$ and a linear coefficient of $-1.$

A \textbf{generating function} encodes a sequence of numbers as coefficients of a formal power series, allowing for manipulation and analysis of the sequence as a single entity \cite{Stanley2015}.
It's crucial in combinatorics for solving problems involving sequences, counting objects, and finding closed-form formulas. 

Typically a generating function would be of the form $\sum_{k\ge 0} a_k t^k$ to encode the sequence $[a_0, a_1, \ldots]$.
In WR, we have more than a sequence, we have an array of potentially infinite dimensions \cite{Wallis1685}. In that case the generating function takes the form $\sum_{m_2, m_3, \ldots \ge 0 \ \ } C[m_2, m_3, \dots] t_2^{m_2} t_3^{m_3} \cdots $ where each array value is the coefficient of a monomial $t_2^{m_2} t_3^{m_3} \cdots $ that encodes its index.

In WR, the particular array of interest is called the \textbf{hyper-Catalan numbers}.
The hyper-Catalan number $C[m_2,m_3,m_4,\ldots]$ counts the ways to subdivide a planar polygon into exactly $m_2$ triangles, $m_3$ quadrilaterals, $m_4$ pentagons, etc.
This is a generalization of the Catalan story 
\cite{Pak2015}, which began in 1751 when Euler \cite{Euler1751} counted triangulations of a polygon; the number of ways one can split a polygon into triangles by drawing diagonals. I've performed a thorough literature review of the Catalan numbers to gain a better understanding of this story.
\cite{Graham1989}

\section{\fontsize{14}{16}\selectfont \textbf{Subdigons}}

We define a \textbf{subdigon} as a roofed polygon subdivided by non-crossing diagonals. 
The \textbf{roof} is a distinguished edge that we draw in red.
The unsubdivided polygon containing the roof is called the \textbf{central polygon}.
As mentioned before, the hyper-Catalan number counts the number of subdigons that has $m_2$ triangles, $m_3$ quadrilaterals, etc \cite{Roman2015}. The type $\m=[m_2, m_3, m_4, \ldots]$ of a subdigon counts the number of triangles, quadrilaterals, etc. it has. Additionally, we define $|$ as the null subdigon, with the $\m=[\ ]$; all $m_k=0$.

We're sorting an infinite collection of subdigons into piles with a finite number of items, where we group all subdigons of type $\m$ into one pile. Thus, the generating function is all subdigons grouped into its unique type (each type represented by monomial $t_2^{m_2} t_3^{m_3} t_4^{m_4} \cdots$), with a finite number of subdigons in each type \cite{Erdelyi1940}. The generating function is made up of the sum of all these monomials.  The coefficient on a given monomial counts the number of subdigons whose type is represented by the monomial.

\begin{figure}[H]
\centering
\includegraphics[width=1.0\textwidth]{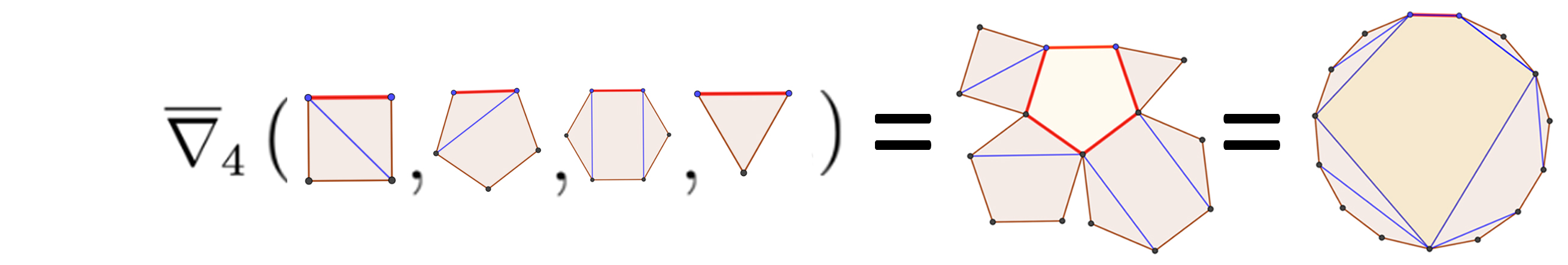}  %
\vspace{-20pt}
\caption{The algebra of subdigons: $\tri_4$  creates a subdigon with a central pentagon.}
\label{fig:tri4} 
\end{figure}

The set of subdigons supports a family of $k$-ary operators. We now define $\tri_k(s_1, s_2, \ldots, s_k)$ as the operation that merges $k$ subdigons together, where $k$ represents a ($k+1$)-sided central polygon. For each side of this polygon (excluding the roof), we will place a subdigon such that its roof coincides with its corresponding side. Note that the subdigons are attached in counter-clockwise order (Figure \ref{fig:tri4}).

However, can we reverse this operation? That is, does there exist a unique set of subdigons which you can apply this $\tri$ operation to to create any subdigon? First, we look at the polygon that has the roof of the subdigon as one of its sides. This is our central polygon; let's say it has $k$ sides. We can then use the $\tri_{k-1}$ operator, with parameters that are the subdigons that are attached to the central polygon, starting from the left of the roof and going in counterclockwise order. We can uniquely create any subdigon using this process as the central polygon will always be unique, which defines what the $k-1$ subdigons will be. This is the basis of constructing our proof.

Now we can define $\Sm_\m$ as the multiset of subdigons of type $\m$; it is necessarily finite, and its size is \textbf{the hyper-Catalan number} $\mC=C[m_2, m_3, \ldots]$.  Then we can define the infinite multiset of subdigons as:
$$\Sm = \sum_{\m \ge \mathbf 0} \Sm_\m .$$

Note that $\m$ is defined as the vector $[m_2, m_3, \dots]$. Let's define the $\psi(s)$ operator as an algebraic monomial that represents the type $\m$ of the subdigon $s$ in the exponents. The $\psi$ operator applied to a subdigon with $m_2$ triangles, $m_3$ quadrilaterals, $m_4$ pentagons, etc. returns: $$\psi(s)=t_2^{m_2}t_3^{m_3}t_4^{m_4}\cdots = \T^{\m}$$
where we use $\T^{\m}$ as an abbreviation for the monomial representing type $\m$.

In general, since $\tri_k$ adds a $(k+1)$-gon we have the identity:
 $$\psi(\tri_k(s_1, s_2, \ldots, s_k))=t_k \, \psi(s_1) \psi(s_2) \cdots \psi(s_k).$$ 

 The multiplications of $\psi(s_i)$ are essentially additions of the type vectors, happening in the exponents.  The outer $t_k$ is also a vector addition, adding $1$ to the $k$th component because the operator $\tri_k$ adds a central $k+1$-gon.

 We extend $\psi$ to multisets of subdigons.  We could use sets, as their elements are never repeated, but multisets are convenient because they have an addition operator: combine the multisets together.  We simply apply $\psi$ to every element and add up the resulting monomials to get a polynomial:
$$\psi(M) = \sum_{s \in M} \psi(s)$$

We can apply $\psi$ to $\Sm_\m$; we get $\psi(\Sm_\m) = \mC \T^\m. $
Thus when we apply $\psi$ to $\Sm$, we get the generating function for the hyper-Catalan array:
$$ \Sb = \Sb[t_2, t_3, \ldots] = \psi(\Sm) = \sum_{\m \ge \mathbf 0} \psi(\Sm_\m) = \sum_{\m \ge \mathbf 0} \mC  \T^{\m}  
 $$

 Note that in $s=\tri_k(s_1, s_2, \ldots, s_k)$, the $s_i$ can be any subdigon; and since every subdigon can be uniquely expressed using our $\tri_k$ operator (or is the null subdigon, a single side), we can create the following multiset equation, essentially a grammar for subdigons:

 $$\Sm =  \ \OB \  | \ \CB \ +  \tri_2(\Sm,\Sm) +  \tri_3(\Sm,\Sm,\Sm)  +  \tri_4(\Sm,\Sm,\Sm,\Sm)  + \ldots$$

 This shows that $\Sm$ can have a central triangle, a central quadrilateral, etc, and its parameters can be any subdigon in S. 
 From our uniqueness argument above, each subdigon appears exactly once on the right.
 
 Let's apply $\psi$ to the equation:

 $$\psi(\Sm) =  \psi(\ \OB \  | \ \CB \ )+  \psi(\tri_2(\Sm,\Sm)) +  \psi(\tri_3(\Sm,\Sm,\Sm))  +  \psi(\tri_4(\Sm,\Sm,\Sm,\Sm))  + \ldots$$

 Simplifying,
 \begin{align}
 \psi(\Sm) &= 1 + t_2\psi(\Sm)^2+t_3\psi(\Sm)^3+t_4\psi(\Sm)^4+\dots
 \\
\Sb &= 1 + t_2 \Sb^2 + t_3 \Sb^3 + t_4 \Sb^4 + \ldots
  \end{align}
We see get a series solution for any polynomial with constant coefficient $1$ and linear coefficient $-1$, namely $\Sb$.  We've proven Theorem \ref{thm:geompoly}.
    
\begin{align}
0=1 - \alpha + t_2  \alpha ^2 + t_3  \alpha^3 + t_4 \alpha^4 + t_5 \alpha^5 \ + \  \ldots 
\end{align}
has a solution
$$\alpha=\psi(\Sm)=\sum_{\substack{m_k \ge 0}} C[m_2, m_3, \dots]   t_2^{m_2 } t_3^{m_3}\dots$$
$$\alpha = \psi(\mathbf{S}) = \sum_{m_k \geq 0} C[m_2, m_3, \ldots] \, t_2^{m_2} t_3^{m_3} \ldots$$

where $C[m_2, m_3, m_4, \dots]$ counts the number of subdigons with $m_2$ triangles, $m_3$ quadrilaterals, $m_4$ pentagons, etc.

That completes our review of WR's proof.  
Note the solution is derived independently from the counting formula for $\mC$.

\section{\fontsize{14}{16}\selectfont \textbf{Hypothesis: This infinite result can be represented as three families of finite results}}

Instead of considering this result as an infinite one, which is unbounded and must be interpreted as a possibly converging function, we truncate the solution to produce something with a finite number of terms. That is, we bound the number of vertices, edges, or faces (omit anything with `degree' greater than $d$) to create a finite expression. However, for faces in particular, bounding the degree of $f$ will still result in an infinite expression.
We can bound the expression further by only considering polynomials up to degree $q$.  Those polynomials only require subdivision into triangles, quadrilaterals, etc. up to $q+1$-gons.  Once we do this, there are only a finite number of types $\m$ at each face level.

\section{\fontsize{14}{16}\selectfont \textbf{Theorem and Proof}}

We want to interpret WR's formula as a finite result at each vertex, edge or bounded face level.
To do so, we add layering variables that explicitly indicate the number of vertices, edges, and faces.

We create an equation that explicitly has the number of vertices, edges, and faces for each existing monomial, or type $\m$. We will do this by making a substitution for the $t_k$ in equation (1). The key insight is that for every triangle we add to an existing subdigon, it adds $1$ vertex (two vertices are shared), $2$ edges (one edge is shared), and $1$ face. In general, for a $k$-sided polygon, it adds $k-2$ vertices, $k-1$ edges, and $1$ face. 

We can calculate the vertices, edges, and faces of a subdigon from its type.
Clearly the number of faces is the sum of the number of triangles, quadrilaterals, etc.  
For vertices and edges, we note gluing on a triangle adds one vertex and two edges, a quadrilateral adds two vertices and three edges, etc.  We have:
\begin{align} \label{eqn:VEF} 
\mF &= m_2 + m_3 + m_4 + \ldots 
\\
\mV &=  2+ m_2 +  2 m_3 + 3 m_4 + \ldots 
\\
\mE &= 1 + 2 m_2 +  3m_3 + 4m_4 + \ldots
\end{align}

\begin{figure} [H]
    \centering
    \includegraphics[width=1\linewidth]{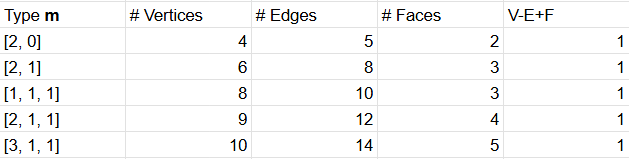}
    \caption{Illustrating Euler's Characteristic Formula in Subdigons}
    \label{fig:Euler}
\end{figure}

We see the Euler Characteristic formula $ V - E + F = 1$ appearing in Figure \ref{fig:Euler}.

Consider the equation $$h(x) = 1 - x + \sum_{k\ge2} v^{k-1} e^k f t_k x^k$$, where every term $x^k$ explicitly has the number of faces, edges, and vertices as powers of $f$, $e$, and $v$, respectively. We want to prove that $h(\alpha)=0$ has a solution with these layered variables. We write $\Sb_L$, our solution layered with vertex, edge, and face variables. 
Let 
\begin{equation}
\Sb_L =\sum_{\substack{m_k \ge 0}} C[m_2, m_3, \dots]   (t_2ve^2f)^{m_2} (t_3v^2e^3f)^{m_3}\dots
\end{equation}
$$\mathbf{S}_{\mathbf{L}} = \sum_{m_k \geq 0} C[m_2, m_3, \ldots] \, (t_2 v e^2 f)^{m_2} (t_3 v^2 e^3 f)^{m_3} \ldots$$

where we replace each $t_i$ with $t_i v^{i-1}e^i f$.
We do this because for each $t_i$ in the expression represents a $i+1$-gon, which adds $i-1$ vertices, $i$ edges, and 1 face.
Expanding this by combining the degrees of $v$, $e$, $f$:
\begin{align}
\Sb_L = &\sum_{m_k \ge 0} C[m_2, m_3, m_4, \dots] t_2^{m_2}t_3^{m_3}t_4^{m_4} \cdots 
\\ & \nonumber \qquad
 v^{m_2+2m_3+3m_4+\dots}
 e^{2m_2+3m_3+4m_4+\dots}
 f^{m_2+m_3+m_4+\dots} 
\\ =& \sum_{\m \ge 0 } \mC \T^{\m} v^{\mV -2}e^{\mE-1}f^{\mF}
\end{align}

We proceed in three steps.  The first is the idea from modular arithmetic that if we evaluate a polynomial on an input and take a remainder at the end, we get the same result if we started by evaluating the polynomial on the remainder and take the remainder at the end.
The second step is to show $h(\Sb_L)=0$, where $h(\alpha)=0$ is the equation we're trying to solve with layering variables added.
The third step uses these to show that our finite equation is satisfied at every vertex, edge, and face level.

\subsection{\fontsize{12}{14}\selectfont Step one}

\begin{theorem}(Modular arithmetic of polynomials)
Let $f(x) = \sum_{k \ge 0} a_k x^k$.
If $a \equiv b \pmod {z^{d+1}}$ then  $f(a) \equiv f(b) \pmod{z^{d+1}}$.
\end{theorem}

We're using $z$ as a generic variable that might be $v, e,$ or $f$. This is a well-known result from modular arithmetic, but for completeness, we include a proof:

\begin{proof}
\begin{align}
   a &\equiv b \pmod {z^{d+1}}
\\
 a_k a^k &\equiv a_k b^k \pmod {z^{d+1}}
\\
\sum_{k\ge 0}  a_k a^k &\equiv  \sum_{k\ge 0} a_k b^k \pmod {z^{d+1}}
\\
f(a) &\equiv f(b) \pmod{z^{d+1}}
\end{align}
\end{proof}
We switch from Gauss's notation to Knuth's \cite{Knuth1997} and write this:
\begin{equation}
     f(a \mod z^{d+1}) \mod z^{d+1} \equiv f(a) \mod z^{d+1}
\end{equation}

\subsection{\fontsize{12}{14}\selectfont Step two}
\begin{theorem}(Layered geometric polynomial solution)
Given our layered polynomial equation 
$$h(x) = 1 - x + \sum_{k\ge2} v^{k-1} e^k f t_k x^k$$
and our layered series $\Sb_L$ defined as above, we have $h(\Sb_L)=0$.
\end{theorem}

\begin{proof}
    Define $c_m$ as the coefficient of $x^m$ in $h(x)$ and $\C^\m$ as $c_2^{m_2}c_3^{m_3}\cdots$. Then, for every $\m$:
\begin{align} \nonumber
\C^{\m}&=c_2^{m_2}c_3^{m_3}\cdots=(t_2ve^2f)^{m_2}(t_3v^2e^3f)^{m_3}\cdots
\\
&=v^{\sum_{k\ge2}{(k-1)m_k}}e^{\sum_{k\ge2}{km_k}}f^{\sum_{k\ge2}{m_k}}\T^\m=v^{\mV-2}e^{\mE-1}f^{\mF}\T^{\m}
\end{align}
    Comparing this to $\Sb_L$, notice that $\Sb_L=\sum_{\m \ge 0 } \mC \C^\m$. Thus, from the soft geometric polynomial formula (Theorem \ref{thm:geompoly}), we have 
    $$h(\Sb_L)=1 - \Sb_L + \sum_{k\ge2} c_k\Sb_L^k=0$$
\end{proof}

\subsection{\fontsize{12}{14}\selectfont Step three}

\begin{theorem}(Finite interpretation of series zeros)
We would like to prove that at every vertex, edge, or face layering of $\Sb_L$, $h(\Sb_L)=0$. In terms of modular arithmetic:
$$
h(\Sb_L \mod v^{d+1}) = 0, \quad
$$
$$
h(\Sb_L \mod e^{d+1}) = 0, \quad
$$
$$
h(\Sb_L \mod f^{d+1}) = 0. \quad
$$
\end{theorem}
\begin{proof}
From Theorem 2, we have:
    $$h(\Sb_L \mod v^{d+1}) \mod {v^{d+1}}\equiv h(\Sb_L) \mod {v^{d+1}}\equiv 0 $$
    $$h(\Sb_L \mod e^{d+1}) \mod {e^{d+1}}\equiv h(\Sb_L) \mod {e^{d+1}}\equiv 0 $$
    $$h(\Sb_L \mod f^{d+1}) \mod {f^{d+1}}\equiv h(\Sb_L) \mod {f^{d+1}}\equiv 0 $$
\end{proof}

As mentioned above, $\Sb_L \mod f^{d+1}$ is an infinite expression.  We further restrict ourselves to polynomials of degree $q$ to make a finite expression.

\begin{theorem}(Finite interpretation of series zeros of bounded face layers)
Given our layered equation of bounded degree,
$$
h_q(x) = 1 - x + \sum_{2 \le k\le q} v^{k-1} e^k f t_k x^k
$$
and our layered solution bounded by degree:
$$
\Sb_{Lq} =\sum_{\substack{m_k \ge 0}} C[m_2, m_3, \ldots, m_q]   (t_2ve^2f)^{m_2} (t_3v^2e^3f)^{m_3}\cdots (t_q v^{q-1} e^q f)^{m_q}
$$
we have
$$
h(\Sb_{Lq} \mod f^{d+1}) = 0 \quad
$$
\end{theorem}
\begin{proof}
We only need up to $q+1$-gon terms to solve polynomials of degree $q$, so we have:
    $$h(\Sb_{Lq} \mod f^{d+1}) \mod {f^{d+1}}\equiv h(\Sb_{Lq}) \mod {f^{d+1}}\equiv 0 $$
\end{proof}

\section{\fontsize{14}{16}\selectfont \textbf{Visualization}}
Next, we aim to better understand how the $\tri$ operator functions through animations. To achieve this, We implemented Python \cite{Sympy2017} code using Jupyter Notebook and Matplotlib that takes any subdigon as input and generates an animation. Since we previously proved that every subdigon can be uniquely expressed as $\tri_k(s_1, s_2, \ldots, s_k)$, the animation visualizes how $\tri_k$ merges $k$ subdigons. The animation consists of four main stages:
\begin{enumerate}
\item The first frame of the animation displays the unique set of $(s_1, s_2, \ldots)$ such that $\tri_k(s_1, s_2, \ldots, s_k)$ is the inputted subdigon. See Figure \ref{fig:stage1} below.
\item Morphing the $\tri_k$ operator into a ($k+1$)-sided central polygon above the $k$ subdigons. We are now ready to merge the subdigons. See Figure \ref{fig:stage2} below.
\item Merging the $k$ subdigons to the ($k+1$) central polygon. See Figure \ref{fig:stage3} below.
\item Morphing the vertices of this combined structure into those of a regular polygon. See Figure \ref{fig:stage4} below.
\end{enumerate}
\begin{figure} [H]
    \centering
    \includegraphics[width=0.5\linewidth]{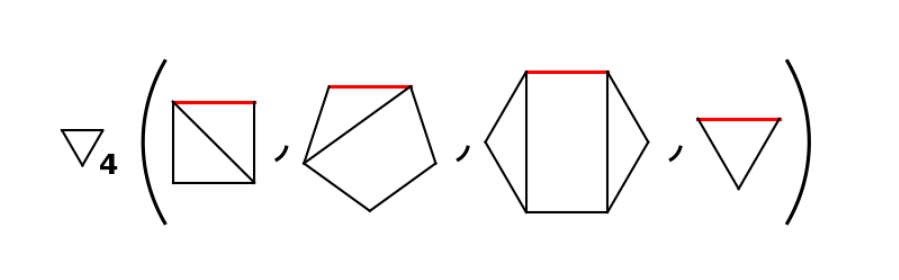}
    \caption{Stage 1: First Frame of the Animation}
    \label{fig:stage1}
\end{figure}
\begin{figure} [H]
    \centering
    \includegraphics[width=0.5\linewidth]{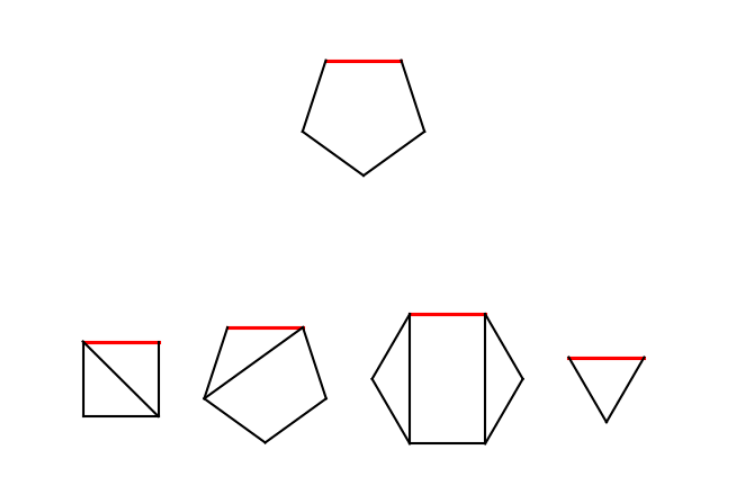}
    \caption{Stage 2: About to Merge Subdigons; Frame 20}
    \label{fig:stage2}
\end{figure}
\begin{figure} [H]
    \centering
    \includegraphics[width=0.5\linewidth]{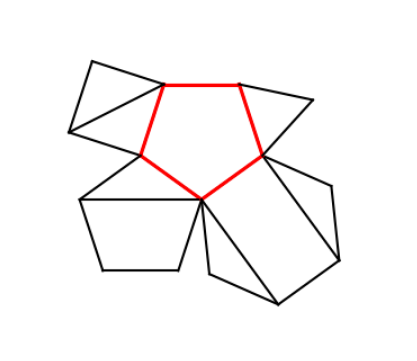}
    \caption{Stage 3: Subdigons Merged; Frame 50}
    \label{fig:stage3}
\end{figure}
\begin{figure} [H]
    \centering
    \includegraphics[width=0.5\linewidth]{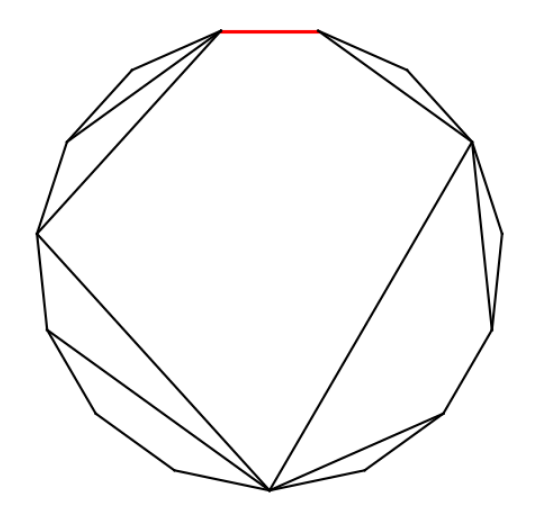}
    \caption{Stage 4: Shifted Vertices; Final Product; Frame 75}
    \label{fig:stage4}
\end{figure}
To create a smooth animation, we used Matplotlib to move objects using a weighted average method. For each object's movement, we defined arrays of starting and ending coordinates between stages. The current coordinates were calculated based on these arrays, where the position is a fraction of the starting coordinates plus a fraction of the ending coordinates. As the frame number increases, we "weigh" the starting array less and the ending array more. For example, in a 21-frame animation:
\begin{itemize}
    \item $i=0: c=1*s+0*e$
    \item $i=6: c=0.7*s+0.3*e$
    \item $i=7: c=0.65*s+0.35*e$
    \item $i=13: c=0.35*s+0.65*e$
    \item $i=19: c=0.05*s+0.95*e$
    \item $i=20: c=0*s+1*e$
\end{itemize}
Here, $c$, $s$, and $e$ represent the current, starting, and ending position arrays, respectively. In the initial frames, the position is more influenced by $s$, while in the later frames, it is increasingly influenced by $e$. This gradual shift in weight creates the illusion of smooth object movement.
\newpage 

\section{\fontsize{14}{16}\selectfont \textbf{Applications}}
The series solution to polynomial equations has several real world impacts. In engineering and control systems, behavior is often modeled with differential equations that reduce to polynomial equations. Having a series solution to these equations (or at least an approximation) allow engineers to design controllers with an optimized performance when designing and analyzing systems like circuits, mechanical structures, or control systems.
Many quantum mechanics problems are approached via perturbation theory, which is a method that relies on a series solution. For example, energy levels and wavefunctions are approximated using series expansions when exact solutions are not feasible to create. Series solutions also help in understanding how small changes in a system can affect its overall behavior.

\section{\fontsize{14}{16}\selectfont \textbf{Future Work + Conclusion}}
There are still many unknown properties about $\Sb$, as well as this derived expression $\Sb_L$. These polygons are nth degree polynomials and should have n complex roots. Our work here only describes one (real) root, but how would we derive a series solution for any/all of them? Is it even possible to write a series for any of the n roots? How can we use this for numerical approximations? When does this expression converge? Computers cannot compute an infinite expression; eventually they are going to have to cut it short, and our work proves that regardless at which vertex, edge, or face you cut it short at, this truncated expression is still a root of the polynomial up to that degree.

We've shown that our variant of the infinite result created by Wildberger and Rubine that explicitly has the number of vertices, edges, and faces in each monomial can be interpreted as three families, each with finite levels through truncating the expression $\Sb_L$. Additionally, we created visualizations depicting why the geometric polynomial formula works by illustrating how the $\tri_k$ operator combines $k$ subdigons and a ($k+1$)-sided central polygon for $k \ge 2$.

\newpage
\bibliographystyle{vancouver}
\bibliography{HyperCatBib.bib}

\end{document}